\documentclass[10pt]{llncs}
\usepackage{amssymb}
\usepackage{latexsym}
\usepackage{amsfonts}
\usepackage{amsmath}
\usepackage{amssymb}

\newcommand{\SL}[2]{{\sf SL}(#1,#2)}
\newcommand{\Sp}[2]{{\sf Sp}(#1,#2)}

\newcommand{\SU}[2]{{\sf SU}(#1,#2)}

\newcommand{\X}{\mathcal X}
\newcommand{\F}{\mathbb F}
\newcommand{\conj}{\,\widehat{\ }\,}
\begin{document}

\title{Constructive membership testing in black-box classical groups}
\author{Sophie Ambrose\inst{1}, Scott H. Murray\inst{2}, Cheryl E. Praeger\inst{1}, and 
Csaba Schneider\inst{3}}
\institute{
School of Mathematics and Statistics, The University of Western Australia \\
35 Stirling Highway
CRAWLEY WA 6009 Australia\\
\email{Cheryl.Praeger@uwa.edu.au,\ alias.sqbr@gmail.com}
\and
Faculty of Information Sciences and Engineering, The University of Canberra, ACT, 2601\\
\email{Scott.Murray@canberra.edu.au}
\and 
Centro de \'Algebra da Universidade de Lisboa\\
Av. Prof. Gama Pinto, 2, 1649-003 Lisboa, Portugal\\
\email{csaba.schneider@gmail.com}}

\maketitle



The research described in this note aims at solving the 
constructive membership problem for the class of quasisimple
classical groups. Our algorithms are developed in 
the black-box group model (see~\cite[Section~3.1.4]{hcgt}); that is,
they do not require specific characteristics of the representations in
which the input groups are given. The elements
of a black-box group are represented, not necessarily uniquely, 
as bit strings of uniform length. We assume the existence of 
oracles to compute the product of two elements, 
the inverse of an element, and to test if two strings represent
the same  element. 
Solving the {\em constructive membership
problem} for a black-box group $G$ requires to write every element
of $G$ as a word in a given generating set. In practice we write the elements
of $G$ as straight-line programs (SLPs) which can be viewed as a compact way
of writing words; see~\cite[Section~3.1.3]{hcgt}.

The constructive membership problem is one of the main tasks identified
in the matrix group recognition project; see~\cite{ob2} for 
details.

The goal of our research is to develop and implement algorithms
to solve the constructive
recognition problem in the classes of black-box classical groups. 
The same problem was already treated by~\cite{ks}.
The main difference between our approach and that 
of~\cite{ks} is that we use the standard generating set of classical groups
given in~\cite{obclg} instead of the larger generating set in~\cite{ks}
and that our goal is to develop algorithms that, in addition to having
good theoretical complexity, perform
well in  practice. Another related algorithm is that of Costi's~\cite{costi}
that solves the constructive membership problem for matrix representations
of classical groups in the defining characteristic. In our algorithms we
reduce the more general problem to a case treated by Costi; see Step~4 below.

In order to briefly explain the main steps of our procedures, 
we use $\Sp {2n}q$ as an example. The natural copy of $\Sp{2n}q$ is the 
group of $2n\times 2n$ matrices over $\F_q$ that preserve a given 
(non-degenerate) symplectic form of a vector space $V=V(\F_q,2n)$. The elements 
of $\Sp{2n}q$ are considered with respect to a given basis 
$e_1,\ldots,e_n,f_n,\ldots,f_1$ consisting of hyperbolic pairs $(e_i,f_i)$. 
The standard generating set $\{\overline s, \overline t, \overline \delta,
\overline u,\overline v,\overline x\}$ of $\Sp{2n}q$ with odd $q$ 
is described in \cite{obclg}. Let $\omega$ be a fixed primitive element
of $\F_q$. Then the standard generators are as follows:
$\overline s:e_1\mapsto f_1,\ f_1\mapsto -e_1$; $\overline t: e_1\mapsto e_1+f_1$; $\overline \delta: e_1\mapsto \omega e_1,\ f_1\mapsto \omega^{-1}f_1$; 
$\overline u: e_1\leftrightarrow e_2,\ f_1\leftrightarrow f_2$; 
$\overline v:e_1\mapsto e_2\mapsto\cdots\mapsto e_n\mapsto e_1,\ f_1\mapsto f_2\mapsto\cdots\mapsto f_n\mapsto f_1$;
$\overline x: f_1\mapsto e_1+f_1,\ f_2\mapsto f_2+e_2$. The standard generators
fix the basis vectors whose images are not listed.

Suppose that $G$ is a black-box group
that is known to be isomorphic to $\Sp{2n}q$ with given $n$ and $q$ and
let us further assume that a generating set $\X=\{s,t,\delta,u,v,x\}$  is
identified in $G$ such that the map $\overline s\mapsto s$, 
$\overline t\mapsto t$, $\overline \delta \mapsto \delta $, 
$\overline u\mapsto u$, $\overline v\mapsto v$, $\overline x\mapsto x$
extends to an isomorphism $\Sp{2n}q\rightarrow G$. This can be achieved using
the algorithms described in~\cite{obclg}.
Our aim is to write a given element $g$ of $G$ 
as an SLP in $\X$. For $g\in G$ let $\widetilde g$ denote the preimage
of $g$ in $\Sp {2n}q$ and let $\widetilde g_{i,j}$ denote the
$(i,j)$-entry of the matrix $\widetilde g$. In order to avoid conjugate 
towers the element $a^b=b^{-1}ab$ will be denoted by $a\conj b$.
Suppose that $q$ is odd, set $\F=\F_q$. 
Our procedure is split into several steps.

{\bf Step 1.} Set $S=\{g\in G\ |\ \widetilde g_{1,2n}=0\}$. 
In this step we find an element $z\in G$ as an SLP in $\X$ 
such that $gz\in S$. 
Set $q=t\conj s$.
For $h\in G$, we have that  $h\in S$ if and only if $q^{q^h}=q$,
and thus we obtain a black-box membership test in $S$ using $O(1)$ 
black-box operations. Since the elements $s$, $u$, and $v$ induce a
transitive group on the subspaces $\left<e_i\right>$, $\left<f_i\right>$
this test can be used to test if $\widetilde g_{1,i}=0$ for all 
$i\in\{1,\ldots,2n\}$. 
If $g\in S$ then we can choose $z=1$; hence assume  that $g\not\in S$.
For $\alpha\in\F$ let $z_\alpha$ be the element of $G$ that 
corresponds to the transformation that maps $e_1\mapsto 
e_1-\alpha e_2$, $f_2\mapsto \alpha f_1+f_2$, and fixes the other 
basis elements.
If 
$\widetilde g_{1,n-1}\neq 0$, 
then $gz_\alpha\in S$ with
$\alpha=-\widetilde g_{1,n}/\widetilde g_{1,n-1}$. Using that 
$z_1=x\conj s$ and $z_{\omega^k}=z_1\conj(\delta^{-k})$, the elements 
$z_{\alpha}$ ($\alpha\in\F$) can be enumerated using $O(q)$ black-box operations and,
for each such $z_{\alpha}$, we can test if $gz_{\alpha}\in S$ using $O(1)$ 
black-box operations. If $gz_{\alpha}\not\in S$ for all $\alpha\in\F$, then
we conclude that $\widetilde g_{1,n-1}=0$, and so $gu\in S$. 
Therefore the cost of finding the suitable $z$ is $O(q)$ black-box
operations. As $z_{\omega^k}=(x\conj s)\conj (\delta^{-k})$, using fast 
exponentiation, the required element $z$ can be written as a SLP of length $O(\log q)$. 

{\bf Step 2.} In this step we assume that $g\in S$ where $S$ is the subset
defined in Step~1. Let $T$ denote the stabilizer 
of the subspace $\left<e_1\right>$. We may assume that $\widetilde 
g_{1,1}\widetilde g_{1,2n-1}\neq 0$ as this can be achieved using the 
membership test in Step~1 with $O(n)$ black-box operations.
We want to find an element $z$ as an SLP in $\X$ 
such that $gz\in T$. 
If 
$(\alpha_2,\ldots,\alpha_{n},\beta_n,\ldots,\beta_1)\in \F^{2n-1}$ 
then let $t(\alpha_2,\ldots,\alpha_n,\beta_n,\ldots,\beta_1)$
denote the element of $G$ corresponding to the transformation
that maps $e_1\mapsto e_1+\alpha_2e_2+\ldots+\alpha_ne_n+\beta_nf_n
+\ldots+\beta_1f_1$,
$e_i\mapsto e_i-\beta_if_1$, $f_i\mapsto f_i+\alpha_if_1$ 
if $i\in\{2,\ldots,n\}$, and $f_1\mapsto f_1$. 
Note that $gt(-\widetilde g_{1,2}/\widetilde g_{1,1},\ldots,-\widetilde g_{1,2n-1}/\widetilde g_{1,1},0)\in T$.
Set $b=(x\conj((t\conj s)\conj g)x^{-1})\conj s$. Then $b=t(\gamma_2,\ldots,\gamma_{2n})$ with
$\gamma_i=-\widetilde g_{1,i}\widetilde g_{1,2n-1}$ for $i=2,\ldots,2n-1$, 
and
$\gamma_{2n}=-\widetilde g_{1,2n-1}(2\widetilde g_{1,1}-\widetilde g_{1,1}^2\widetilde g_{1,2n}-\widetilde g_{1,2n-1})$.
Further, $b\conj (\delta^{-k})=t(\gamma_1\omega^k,\ldots,\gamma_{2n-1}\omega^k,\gamma_{2n}\omega^{2k})$.
Hence there is some $k_0$ such that $\gamma_i\omega^{k_0}=-\widetilde g_{1,i}/\widetilde g_{1,1}$ for all $i\in\{2,\ldots,2n-1\}$. Set $z_0=b\conj(\delta^{-k_0})$. 
Using
the membership test for $S$ explained in the previous paragraph
and using the fact that 
$(\widetilde{gz_0})_{1,2n-1}=0$,
the element $z_0$ can be found using $O(q)$ group multiplications.
Now given the element $z_0$, we can recover the entries $\gamma_i\omega^{k_0}$
and we can write the element $z_0=t(\gamma_2\omega^{k_0},\ldots,\gamma_{2n-1}\omega^{k_0})$ as SLP
as follows. For $i=2,\ldots,n-1$ and $\alpha\in \F$
let $x_i(\alpha)=t(0,\ldots,0,\alpha,0,\ldots,0)$ where the non-zero
entry appears in the $(i-1)$-th position.
Let $I$ denote the set of indices $i\in\{2,\ldots,2n-1\}$ for which
$\gamma_i\neq 0$.  For $i\in I$, let $k_i$ be such
that $\gamma_i\omega^{k_0}=\omega^{-k_i}$. Then 
$z_0=\prod_{i\in I}x_i(\omega^{k_i})$. 
As 
$x_2(1)=(x\conj s)^{-1}$, $x_{i+1}(1)$ can be obtained from $x_{i}(1)$ 
using $O(1)$ group multiplications, and
$x_i(1)\conj(\delta^{-k})= x_i(\omega^k)$, the entries $\gamma_i\omega^{k_0}$ can be
recovered and $z_0$ can be written as an SLP using $O(nq)$ group 
multiplications. The length of the SLP is $O(n\log q)$.

{\bf Step 3.} Now we assume that $g\in T$. We repeat Steps~1 and~2 with 
$g^s$ and obtain an element $z_l$ and $z_r$ as SLPs in $\X$ 
 such that $z_lgz_r$ is in the intersection $G_1$ of the stabilizers of 
$\left<e_1\right>$ and $\left<f_1\right>$. The cost of this step is $O(nq)$
group multiplication and the length of the SLPs to $z_l$ and $z_r$ is 
$O(n\log q)$. 

{\bf Step 4.}
In this step we assume that $g\in G$ lies in $G_1$.
We have, for $i\in\{2,\ldots,2n-1\}$,
that $x_i(1)\conj g=t(\widetilde g_{i,2}/\widetilde g_{1,1},\ldots,\widetilde 
g_{i,2n-1}/g_{1,1})$. Using the procedures described in Step~2, 
the entries $\widetilde g_{i,j}$ with $i,\ j\in\{2,\ldots,2n-1\}$ 
 can be recovered using $O(n^2q)$ multiplications.
Let $M$ denote the $(2n-2)\times (2n-2)$ matrix formed by these entries.
The procedure of Costi~\cite{costi}
is used to write $M$ as
a SLP in the standard generators of $\Sp{2n-2}q$. Considering 
$\Sp{2n-2}q$ as the subgroup $G_1$, we evaluate this SLP in $G$ 
to obtain an element $z$. Now $gz^{-1}$ is a diagonal matrix with
$(\widetilde{gz^{-1}})_{2,2}=\cdots=(\widetilde{gz^{-1}})_{2n-1,2n-1}$. Hence 
there is some $k$ such that $gz^{-1}\delta^k\in Z(G)$. 

After the end of Step~4, the element $g$ is written as an SLP in $\X$ modulo
the center of $G$ using $O(n^2q)$ group operations. The length of the SLP is 
$O(n^2\log q)$. We emphasize that the procedures, 
while using vector and matrix notation for ease of the exposition, 
are actually black-box procedures requiring only the basic group 
operations 
of multiplication, 
inversion and equality testing.
Similar algorithms are developed and implemented for 
the classical groups $\SL nq$  and  $\SU nq$ (with odd $q$). The implementations
are available in the computational algebra system {\sc Magma}.

{\sc Acknowledgment.} 
Murray would like to thank the
Magma project at the University of Sydney, where some of the work was
carried out. 
Praeger would like to acknowledge  the support of the
Australian Research Council Discovery Grant DP0879134 and Federation Fellowship FF0776186. Schneider was supported by the FCT project
PTDC/MAT/101993/2008 (Portugal) and by the OTKA grant~72845 (Hungary).


\begin{thebibliography}{XXXXX}

\bibitem[Cos]{costi}
Elliot Costi. 
\newblock
Constructive membership testing in classical groups. 
\newblock PhD thesis, 
Queen
Mary, University of London, 2009.



\bibitem[HCGT]{hcgt}
Derek F. Holt, Bettina Eick, Eamonn A. O'Brien.
\newblock {\em The Handbook of Computational Group Theory.}
\newblock Chapman and Hall/CRC, 2005.

\bibitem[KS]{ks}
William M. Kantor and \'Akos Seress. {\em Black box classical groups}.
Vol.~149 of {\em Memoirs Amer. Math. Soc.}, 2001.

\bibitem[LGOB]{obclg}
C. R. Leedham-Green and E. A. O'Brien. 
\newblock 
Constructive recognition of classical groups in odd characteristic.
\newblock {\em Journal of Algebra} {\bf 322}, 833-881, 2009.

\bibitem[OB]{ob2}
E. A. O'Brien. 
\newblock
Algorithms for matrix groups.
\newblock Groups St Andrews, (Bath), August 2009, accepted to appear 2010.
\end{thebibliography}
\end{document}